\documentclass[12pt]{article}
\usepackage{latexsym,amsfonts,amsmath,amssymb}
\newtheorem{theorem}{Theorem}
\newtheorem{corollary}[theorem]{Corollary}
\newtheorem{sublemma}{Lemma}[theorem]
\newtheorem{lemma}[theorem]{Lemma}
\newtheorem{question}[theorem]{Question}
\newtheorem{observation}[theorem]{Observation}
\newtheorem{claim}[theorem]{Claim}
\newtheorem{subclaim}{Claim}[sublemma]
\newtheorem{conjecture}[theorem]{Conjecture}
\newtheorem{fact}[theorem]{Fact}
\newtheorem{definition}[theorem]{Definition}
\newtheorem{remark}[theorem]{Remark}
\newtheorem{example}[theorem]{Example}
\newtheorem{exercise}{Exercise}[section]
\def\Theorem #1.#2 #3\par{\setbox1=\hbox{#1}\ifdim\wd1=0pt
   \begin{theorem}{\rm #2} #3\end{theorem}\else
   \newtheorem{#1}[theorem]{#1}\begin{#1}\label{#1}{\rm #2} #3\end{#1}\fi}
\def\Corollary #1.#2 #3\par{\setbox1=\hbox{#1}\ifdim\wd1=0pt
   \begin{corollary}{\rm #2} #3\end{corollary}\else
   \newtheorem{#1}[theorem]{#1}\begin{#1}\label{#1}{\rm #2} #3\end{#1}\fi}
\def\Lemma #1.#2 #3\par{\setbox1=\hbox{#1}\ifdim\wd1=0pt
   \begin{lemma}{\rm #2} #3\end{lemma}\else
   \newtheorem{#1}[theorem]{#1}\begin{#1}\label{#1}{\rm #2} #3\end{#1}\fi}
\def\SubLemma #1.#2 #3\par{\setbox1=\hbox{#1}\ifdim\wd1=0pt
   \begin{sublemma}{\rm #2} #3\end{sublemma}\else
   \newtheorem{#1}{#1}[theorem]\begin{#1}\label{#1}{\rm #2} #3\end{#1}\fi}
\def\Question #1.#2 #3\par{\setbox1=\hbox{#1}\ifdim\wd1=0pt
   \begin{question}{\rm #2} #3\end{question}\else
   \newtheorem{#1}[theorem]{#1}\begin{#1}\label{#1}{\rm #2} #3\end{#1}\fi}
\def\Observation #1.#2 #3\par{\setbox1=\hbox{#1}\ifdim\wd1=0pt
   \begin{observation}{\rm #2} #3\end{observation}\else
   \newtheorem{#1}[theorem]{#1}\begin{#1}\label{#1}{\rm #2} #3\end{#1}\fi}
\def\Claim #1.#2 #3\par{\setbox1=\hbox{#1}\ifdim\wd1=0pt
   \begin{claim}{\rm #2} #3\end{claim}\else
   \newtheorem{#1}[theorem]{#1}\begin{#1}\label{#1}{\rm #2} #3\end{#1}\fi}
\def\SubClaim #1.#2 #3\par{\setbox1=\hbox{#1}\ifdim\wd1=0pt
   \begin{subclaim}{\rm #2} #3\end{subclaim}\else
   \newtheorem{#1}{#1}[sublemma]\begin{#1}\label{#1}{\rm #2} #3\end{#1}\fi}
\def\Conjecture #1.#2 #3\par{\setbox1=\hbox{#1}\ifdim\wd1=0pt
   \begin{conjecture}{\rm #2} #3\end{conjecture}\else
   \newtheorem{#1}[theorem]{#1}\begin{#1}\label{#1}{\rm #2} #3\end{#1}\fi}
\def\Fact #1.#2 #3\par{\setbox1=\hbox{#1}\ifdim\wd1=0pt
   \begin{fact}{\rm #2} #3\end{fact}\else
   \newtheorem{#1}[theorem]{#1}\begin{#1}\label{#1}{\rm #2} #3\end{#1}\fi}
\def\Definition #1.#2 #3\par{\setbox1=\hbox{#1}\ifdim\wd1=0pt
   \begin{definition}{\rm #2} {\rm #3}\end{definition}\else
   \newtheorem{#1}[theorem]{#1}\begin{#1}\label{#1}{\rm #2} {\rm #3}\end{#1}\fi}
\def\Remark #1.#2 #3\par{\setbox1=\hbox{#1}\ifdim\wd1=0pt
   \begin{remark}{\rm #2} {\rm #3}\end{remark}\else
   \newtheorem{#1}[theorem]{#1}\begin{#1}\label{#1}{\rm #2} {\rm #3}\end{#1}\fi}
\def\Example #1.#2 #3\par{\setbox1=\hbox{#1}\ifdim\wd1=0pt
   \begin{example}{\rm #2} #3\end{example}\else
   \newtheorem{#1}[theorem]{#1}\begin{#1}\label{#1}{\rm #2} #3\end{#1}\fi}
\def\Exercise #1.#2 #3\par{\setbox1=\hbox{#1}\ifdim\wd1=0pt
   {\footnotesize\begin{exercise}{\rm #2} {\rm #3}\end{exercise}}\else
   \newtheorem{#1}[section]{#1}{\footnotesize\begin{#1}\label{#1}{\rm #2} {\rm #3}\end{#1}}\fi}
\def\QuietTheorem #1.#2 #3\par{\setbox1=\hbox{#1}\ifdim\wd1=0pt\proclaim{Theorem {\rm #2}}{#3}\else\proclaim{#1 {\rm #2}}{#3}\fi}
\newcommand{\proclaim}[2]{\smallskip\noindent{\bf #1} {\sl#2}\par\smallskip}
\def\Proclaim #1.#2 #3\par{\proclaim{#1 {\rm #2}}{#3}}
\newenvironment{proof}{\noindent}{\kern2pt\QEDbox\par\bigskip}
\def\Proof#1: {\setbox1=\hbox{#1}\ifdim\wd1=0pt\begin{proof}{\bf Proof: }\else\medskip\begin{proof}{\bf #1: }\fi}
\newcommand{\QED}{\end{proof}}
\def\BF#1.{{\bf #1.}}
\def\Abstract #1\par{\begin{quotation}{\singlespaced\footnotesize{\noindent{\bf Abstract.~}#1}}\end{quotation}}
\def\Title #1\par{\title{#1}\maketitle}
\def\Author #1\par{\author{#1}}
\def\Acknowledgement#1\par{\thanks{#1}}
\def\Chapter #1\par{\chapter{#1}}
\def\Section #1\par{\section{#1}}
\def\QuietSection #1\par{\section*{#1}}
\def\SubSection #1\par{\subsection{#1}}
\def\SubSubSection #1\par{\subsubsection{#1}}
\def\MidTitle #1\par{\bigskip\goodbreak\centerline{\small\bf #1}\bigskip\noindent}
\def\Margin #1\par{\marginpar{\tiny #1}}

\newcommand{\singlespaced}{\baselineskip=15pt}
\def\bottomnote #1\par{{\renewcommand{\thefootnote}{}\footnotetext{#1}}}
\newcommand{\F}{{\mathbb F}}

\renewcommand{\P}{{\mathbb P}}
\newcommand{\Q}{{\mathbb Q}}

\newcommand{\Gtail}{{G_{\!\scriptscriptstyle\rm tail}}}

\newcommand{\Ptail}{{\P_{\!\scriptscriptstyle\rm tail}}}

\newcommand{\Mbar}{{\overline{M}}}
\newcommand{\Nbar}{{\overline{N}}}

\newcommand{\Adot}{{\dot A}}
\newcommand{\Bdot}{{\dot B}}

\newcommand{\one}{\mathop{1\hskip-3pt {\rm l}}}
\newfont{\msam}{msam10 at 12pt}
\newcommand{\from}{\mathbin{\vbox{\baselineskip=3pt\lineskiplimit=0pt
                         \hbox{.}\hbox{.}\hbox{.}}}}
\newcommand{\of}{\subseteq}

\newcommand{\set}[1]{\{\,{#1}\,\}}
\newcommand{\singleton}[1]{\left\{{#1}\right\}}
\newcommand{\compose}{\circ}

\newcommand{\elesub}{\prec}

\newcommand{\muchgt}{\gg}

\newcommand{\dom}{\mathop{\rm dom}}

\newcommand{\ran}{\mathop{\rm ran}}

\newcommand{\add}{\mathop{\rm add}}

\newcommand{\cof}{\mathop{\rm cof}}

\newcommand{\image}{\mathbin{\hbox{\tt\char'42}}}
\newcommand{\plus}{{+}}

\newcommand{\restrict}{\upharpoonright}
\newcommand{\satisfies}{\models}

\newcommand{\cross}{\times}
\newcommand{\concat}{\mathbin{{}^\smallfrown}}

\newcommand{\union}{\cup}

\newcommand{\Union}{\bigcup}

\newcommand{\intersect}{\cap}

\newcommand{\LaverDiamond}{\mathop{\hbox{\line(0,1){10}\line(1,0){8}\line(-4,5){8}}\hskip 1pt}\nolimits}
\newcommand{\LD}{\LaverDiamond}

\newcommand{\LDwc}{\LD^{\hbox{\!\!\tiny wc}}}
\newcommand{\LDunf}{\LD^{\hbox{\!\!\tiny unf}}}
\newcommand{\LDind}{\LD^{\hbox{\!\!\tiny ind}}}

\newcommand{\LDsunf}{\LD^{\hbox{\!\!\tiny sunf}}}

\newcommand{\smalllt}{\mathrel{\mathchoice{\raise2pt\hbox{$\scriptstyle<$}}{\raise1pt\hbox{$\scriptstyle<$}}{\raise0pt\hbox{$\scriptscriptstyle<$}}{\scriptscriptstyle<}}}
\newcommand{\smallleq}{\mathrel{\mathchoice{\raise2pt\hbox{$\scriptstyle\leq$}}{\raise1pt\hbox{$\scriptstyle\leq$}}{\raise1pt\hbox{$\scriptscriptstyle\leq$}}{\scriptscriptstyle\leq}}}

\newcommand{\ltkappa}{{{\smalllt}\kappa}}
\newcommand{\leqkappa}{{{\smallleq}\kappa}}

\newcommand{\leqgamma}{{{\smallleq}\gamma}}

\newcommand{\leqtheta}{{{\smallleq}\theta}}

\newcommand{\leqbeta}{{{\smallleq}\beta}}

\newcommand{\card}[1]{{|#1|}}
\newcommand{\boolval}[1]{\mathopen{\lbrack\!\lbrack}\,#1\,\mathclose{\rbrack\!\rbrack}}
\def\[#1]{\boolval{#1}}

\newcommand{\UnderTilde}[1]{{\setbox1=\hbox{$#1$}\baselineskip=0pt\vtop{\hbox{$#1$}\hbox to\wd1{\hfil$\sim$\hfil}}}{}}
\newcommand{\Undertilde}[1]{{\setbox1=\hbox{$#1$}\baselineskip=0pt\vtop{\hbox{$#1$}\hbox to\wd1{\hfil$\scriptstyle\sim$\hfil}}}{}}
\newcommand{\undertilde}[1]{{\setbox1=\hbox{$#1$}\baselineskip=0pt\vtop{\hbox{$#1$}\hbox to\wd1{\hfil$\scriptscriptstyle\sim$\hfil}}}{}}
\newcommand{\UnderdTilde}[1]{{\setbox1=\hbox{$#1$}\baselineskip=0pt\vtop{\hbox{$#1$}\hbox to\wd1{\hfil$\approx$\hfil}}}{}}
\newcommand{\Underdtilde}[1]{{\setbox1=\hbox{$#1$}\baselineskip=0pt\vtop{\hbox{$#1$}\hbox to\wd1{\hfil\scriptsize$\approx$\hfil}}}{}}

\newcommand{\st}{\mid}
\renewcommand{\th}{{\hbox{\scriptsize th}}}

\newcommand{\iso}{\cong}
\def\<#1>{\langle#1\rangle}

\newcommand{\QEDbox}{\fbox{}}
\newcommand{\cp}{\mathop{\rm cp}}

\newcommand{\REG}{\mathop{\hbox{\sc reg}}}

\newcommand{\INACC}{\mathop{\hbox{\sc inacc}}}

\newcommand{\ZFC}{\hbox{\sc zfc}}

\newcommand{\GCH}{\hbox{\sc gch}}

\newcommand{\inacc}{\hbox{\sc inacc}}

\newcommand{\cell}[1]{\boxit{\hbox to 17pt{\strut\hfil$#1$\hfil}}}
\newcommand{\head}[2]{\lower2pt\vbox{\hbox{\strut\footnotesize\it\hskip3pt#2}\boxit{\cell#1}}}
\newcommand{\boxit}[1]{\setbox4=\hbox{\kern2pt#1\kern2pt}\hbox{\vrule\vbox{\hrule\kern2pt\box4\kern2pt\hrule}\vrule}}
\newcommand{\Col}[3]{\hbox{\vbox{\baselineskip=0pt\parskip=0pt\cell#1\cell#2\cell#3}}}
\newcommand{\tapenames}{\raise 5pt\vbox to .7in{\hbox to .8in{\it\hfill input: \strut}\vfill\hbox to
.8in{\it\hfill scratch: \strut}\vfill\hbox to .8in{\it\hfill output: \strut}}}
\newcommand{\Head}[4]{\lower2pt\vbox{\hbox to25pt{\strut\footnotesize\it\hfill#4\hfill}\boxit{\Col#1#2#3}}}
\newcommand{\Dots}{\raise 5pt\vbox to .7in{\hbox{\ $\cdots$\strut}\vfill\hbox{\ $\cdots$\strut}\vfill\hbox{\
$\cdots$\strut}}}
\renewcommand{\dots}{\raise5pt\hbox{\ $\cdots$}}
\newcommand{\factordiagramup}[6]{$$\begin{array}{ccc}
#1&\raise3pt\vbox{\hbox to60pt{\hfill$\scriptstyle
#2$\hfill}\vskip-6pt\hbox{$\vector(4,0){60}$}}&#3\\ \vbox
to30pt{}&\raise22pt\vtop{\hbox{$\vector(4,-3){60}$}\vskip-22pt\hbox
to60pt{\hfill$\scriptstyle #4\qquad$\hfill}}
     &\ \ \lower22pt\hbox{$\vector(0,3){45}$}\ {\scriptstyle #5}\\
\vbox to15pt{}&&#6\\
\end{array}$$}
\newcommand{\factordiagram}[6]{$$\begin{array}{ccc}
#1&&\\ \ \ \raise22pt\hbox{$\vector(0,-3){45}$}\ {\scriptstyle #2}
&\raise22pt\hbox{$\vector(2,-1){90}$}\raise5pt\llap{$\scriptstyle#3$\qquad\quad}&\vbox
to25pt{}\\ #4&\raise3pt\vbox{\hbox to90pt{\hfill$\scriptstyle
#5$\hfill}\vskip-6pt\hbox{$\vector(4,0){90}$}}&#6\\
\end{array}$$}
\newcommand{\df}{\it} 
\begin{document}
\author{Mirna D\v zamonja\\
\normalsize\sc University of East Anglia\\
\normalsize\sc Norwich, NR4 7TJ, UK\\
\footnotesize http://www.mth.uea.ac.uk/people/md.html\\
\\
Joel David Hamkins\\
\normalsize\sc The City University of New York\\
\normalsize\sc Mathematics, 365 Fifth Avenue, New York, NY 10016\\
\footnotesize http://jdh.hamkins.org\\
}

\bottomnote MSC 03E55, 03E35, 03E05. Keywords: unfoldable cardinals, infinite combinatorics, diamond sequence, Laver function. The research of the
first author has been supported by an EPSRC Advanced Fellowship. The research of the second author has been supported by grants from Georgia State
University, The Research Foundation of The City University of New York and the National Science Foundation. He is affiliated with the College of
Staten Island of CUNY and the CUNY Graduate Center.

\Title Diamond (on the regulars) can fail at any strongly unfoldable cardinal

\Abstract If $\kappa$ is any strongly unfoldable cardinal, then this is preserved in a forcing extension in which $\Diamond_\kappa(\REG)$ fails. This
result continues the progression of the corresponding results for weakly compact cardinals, due to Woodin, and for indescribable cardinals, due to
Hauser.

\Section Introduction

For any stationary subset $E$ of a regular uncountable cardinal $\kappa$ the combinatorial diamond principle $\Diamond_\kappa(E)$ asserts that there
is a sequence $\<A_\alpha\st\alpha<\kappa>$ such that for any $A\of\kappa$ the set $\set{\alpha\in E\st A\intersect\alpha=A_\alpha}$ is stationary in
$\kappa$. It is relatively easy to show, for example, that if $\kappa$ is a measurable cardinal, then $\Diamond_\kappa(\REG)$ holds, where $\REG$ is
the set of regular cardinals (see Observation \ref{MeasurableImpliesDiamond}). The same is true if $\kappa$ is merely ineffable or even subtle, by
results of Jensen and Kunen (see \cite{Kanamori1988:DiamondsLCandUltrafilters}). Furthermore, $\Diamond_\kappa(E)$ holds in $L$ for any stationary
subset $E$ of any uncountable regular cardinal $\kappa$ there. So if a cardinal $\kappa$ (at least weakly Mahlo, so that $\REG$ is stationary)
exhibits too large a large cardinal property, or if the universe is too close to $L$, then one should expect $\Diamond_\kappa(\REG)$ to hold. This
article is about a sequence of results that squeeze between these two requirements, by aiming to produce failures of $\Diamond_\kappa(\REG)$ as high
in the large cardinal hierarchy as possible.

Specifically, our main theorem continues the progression of two earlier results. First, Woodin \cite{CummingsWoodin} showed that
$\Diamond_\kappa(\REG)$ can fail at any weakly compact cardinal; in fact, the existence of a weakly compact cardinal is equiconsistent with the
existence of a weakly compact cardinal $\kappa$ such that $\neg\Diamond_\kappa(\REG)$, plus \GCH. After this, Kai Hauser
\cite{Hauser92:IndescribablesWithoutDiamond} produced failures of $\Diamond_\kappa(\REG)$ at indescribable cardinals, where again the existence of a
$\Pi^m_n$ indescribable cardinal ($m,n\geq 1$) is equiconsistent with the existence of a $\Pi^m_n$ indescribable cardinal $\kappa$ such that
$\neg\Diamond_\kappa(\REG)$, plus \GCH. All these results address a question asked in Kanamori's article
\cite{Kanamori1988:DiamondsLCandUltrafilters}, namely, if it is possible to have a large cardinal $\kappa$ such that $\Diamond_\kappa(\REG)$ fails.
Here, we push the violations of $\Diamond_\kappa(\REG)$ higher with the following:

\QuietTheorem Main Theorem. The existence of a strongly unfoldable cardinal is equiconsistent with the existence of a strongly unfoldable cardinal
$\kappa$ such that $\neg\Diamond_\kappa(\REG)$, plus \GCH.\label{MainEquiconsistency}

In all three cases, what is actually proved is that if a cardinal $\kappa$ is respectively weakly compact, indescribable or strongly unfoldable and
the \GCH\ holds, then this is preserved in a forcing extension in which $\neg\Diamond_\kappa(\REG)$ holds. A firm upper bound on how high in the
large cardinal hierarchy one can have failures of $\Diamond_\kappa(\REG)$ is provided by the following easy observation about measurable cardinals.
(Finer results were obtained by Jensen and Kunen, as we mentioned above.)

\Observation. If $\kappa$ is a measurable cardinal,
 then $\Diamond_\kappa(\REG)$ holds.\label{MeasurableImpliesDiamond}

\Proof: This well known argument follows the usual proof of $\Diamond$ in $L$. If $A_\alpha$ is defined for all $\alpha<\gamma$, then for regular
$\gamma$ let $A_\gamma$ be any subset of $\gamma$ such that $\set{\alpha\in\REG\intersect\gamma\st A_\gamma\intersect\alpha=A_\alpha}$ is not
stationary in $\gamma$, if there is any such set, and use any subset of $\gamma$ otherwise. At every stage of the construction, therefore, we have
added to the sequence $\vec A$ a set which has not yet been sufficiently anticipated, if such a set existed. If the resulting sequence $\vec
A=\<A_\alpha\st\alpha<\kappa>$ does not witness $\Diamond_\kappa(\REG)$, then there is a set $A\of\kappa$ which is not anticipated on a stationary
set of regular cardinals. Let $j:V\to M$ be the ultrapower embedding generated by a normal measure on $\kappa$, and consider the sequence $j(\vec
A)$. Specifically, let $A^*=j(\vec A)(\kappa)$ be the $\kappa^\th$ set in this sequence. By definition, this set should be a subset of $\kappa$ which
is not anticipated by $j(\vec A)\restrict\kappa=\vec A$ on a stationary set of regular cardinals in $M$, if there are any such sets in $M$. But since
$M$ and $V$ have the same subsets of $\kappa$, the two models agree on the regular cardinals up to $\kappa$, on whether a set of such cardinals is
stationary and on whether such a set is anticipated by $\vec A$. And since we assumed that there are subsets of $\kappa$ that are not anticipated by
$\vec A$ on a stationary set in $V$, it follows that there are such sets in $M$. Thus, in particular, the set $A^*=j(\vec A)(\kappa)$ is not
anticipated by $\vec A$ on a stationary set of regular cardinals. So there is a club set $C\of\kappa$ such that if $\gamma\in C\intersect\REG$, then
$A_\gamma\not=A^*\intersect\gamma$. But since $\kappa\in j(C)$, this implies that $j(\vec A)(\kappa)\not=j(A^*)\intersect\kappa$, which is absurd
since $j(A^*)\intersect\kappa=A^*=j(\vec A)(\kappa)$.\QED

If $\kappa$ is measurable, then the diamond sequence $\vec A$ constructed above will witness $\Diamond_\kappa(\REG)$ in $M$, whenever $j:V\to M$ has
critical point $\kappa$, and so the set of cardinals $\gamma<\kappa$ that are weakly compact and satisfy $\Diamond_\gamma(\REG)$ has measure one with
respect to every normal measure on $\kappa$. The construction of Observation \ref{MeasurableImpliesDiamond} has little to do with $\REG$, and it
produces a sequence $\vec A$ whose restriction to any set $E$ in a normal measure on $\kappa$ witnesses $\Diamond_\kappa(E)$.

One can also show that every measurable cardinal $\kappa$ has a kind of Laver function $\ell:\kappa\to V_\kappa$ for weakly compact embeddings,
meaning that for any $A\in H(\kappa^\plus)$ and any transitive set $M$ of size $\kappa$ containing $A$ and $\kappa$, there is a weak compactness
embedding $j:M\to N$ such that $j(\ell)(\kappa)=A$. This is called the Laver diamond principle for weak compactness and denoted $\LDwc_\kappa$ in
\cite{Hamkins:LaverDiamond}, where similar notions are considered for a variety of large cardinals. All of the Laver diamond principles easily imply
$\Diamond_\kappa(\REG)$.

The question of the exact boundary in the large cardinal hierarchy where $\Diamond_\kappa(\REG)$ can fail remains open. The results of this article
squeeze this boundary to somewhere between the strongly unfoldable cardinals and the subtle cardinals.

We use Villaveces' \cite{Villaveces98} embedding characterization of the strongly unfoldable cardinals, through which they resemble miniature strong
cardinals in the same way that embeddings help weakly compact cardinals to resemble miniature measurable cardinals. In each case, an embedding
$j:V\to W$ defined on all of $V$ is replaced with embeddings $j:M\to N$ defined only on a transitive structure of size $\kappa$. In our proof,
however, we also show that strongly unfoldable cardinals have an embedding characterization analogous to the supercompact cardinals, making
``miniature strong'' equivalent to ``miniature supercompact,'' a fact we found surprising.

\Section Strongly Unfoldable Cardinals

The unfoldable and strongly unfoldable cardinals were introduced by Villaveces in \cite{Villaveces98} as a direct generalization of the weakly
compact cardinals. The embedding characterizations will be the most convenient: a cardinal $\kappa$ is {\df weakly compact} if and only if for any
transitive structure $M$ of size $\kappa$ containing $\kappa$ as an element there is a transitive set $N$ and an elementary embedding $j:M\to N$ with
critical point $\kappa$. Since one can always restrict an embedding from a larger domain to a smaller domain, this is equivalent to insisting merely
that every subset $B\of\kappa$ can be placed into a structure $M$ having an appropriate embedding $j:M\to N$. And of course there are such structures
$M$ such that (i) $M$ is transitive of size $\kappa$, (ii) $M$ models some fixed large finite fragment of set theory $\ZFC^*$ and (iii)
$M^\ltkappa\of M$. We will refer to $M$ having these properties as the {\df $\kappa$-models} of set theory. Thus, $\kappa$ is weakly compact if and
only if every $\kappa$-model of set theory $M$ has an embedding $j:M\to N$ into a transitive $N$ with critical point $\kappa$. For the remainder of
this article, when $M$ is transitive and we say that there is an embedding $j:M\to N$, we mean to imply that $N$ is also transitive, and we will
usually also mean that the critical point of $j$ is whatever cardinal $\kappa$ is under consideration. Throughout we shall use the adjective
inaccessible to mean a strongly inaccessible cardinal.

A cardinal $\kappa$ is {\df $\theta$-unfoldable} when every $\kappa$-model of set theory $M$ has an embedding $j:M\to N$ with critical point $\kappa$
and $j(\kappa)\geq\theta$. This is equivalent to insisting that every set $B\of\kappa$ can be placed into a $\kappa$-model $M$ having such an
embedding. We will refer to such embeddings as the $\theta$-unfoldability embeddings. The cardinal $\kappa$ is {\df unfoldable} if it is
$\theta$-unfoldable for every $\theta$. Generalizing this further, $\kappa$ is {\df $\theta$-strongly unfoldable} if every $\kappa$-model of set
theory $M$ has a $\theta$-unfoldability embedding $j:M\to N$ that also satisfies $V_\theta\of N$. And again, this is equivalent to requiring that any
set $B\of\kappa$ can be placed into such an $M$. Finally, of course, $\kappa$ is {\df strongly unfoldable} if it is $\theta$-strongly unfoldable for
every $\theta$. The strongly unfoldable cardinals therefore look a bit like miniature strong cardinals, with embeddings defined only on structures of
size $\kappa$, rather than the entire universe.

If $V=L$, the requirement that $V_\theta\of N$ amounts merely to $L_{\beth_\theta}\of N$, which is true of any $N$ containing enough ordinals, and so
in $L$ every unfoldable cardinal is strongly unfoldable. Since every unfoldable cardinal is unfoldable in $L$ (see \cite{Villaveces98}), the two
unfoldability notions have the same consistency strength. Nevertheless, results in \cite{Hamkins2001:UnfoldableCardinals} show that the strong
unfoldability of any unfoldable cardinal can be destroyed by forcing that preserves its unfoldability, so as large cardinal notions the two concepts
are distinct.

The next two theorems show that the consistency strength of the existence of unfoldable or strongly unfoldable cardinals lies strictly between the
totally indescribable cardinals and the subtle cardinals. Thus, our Main Theorem finds violations of $\Diamond_\kappa(\REG)$ higher in the large
cardinal hierarchy than \cite{Hauser92:IndescribablesWithoutDiamond}.)

By \cite{Hauser1991:IndescribableCardinals}, a cardinal $\kappa$ is {\df $\Pi^m_n$-indescribable} if for any $\kappa$-model $M$ there is an embedding
$j:M\to N$ with critical point $\kappa$ such that $N$ is $\Sigma^m_n$-correct, that is, $(V_{\kappa+m})^N\elesub_n V_{\kappa+m}$ and
$N^\card{V_{\kappa+m-2}}\of N$ (meaning $N^\ltkappa\of N$ when $m=1$). The usual reflection definition shows that it is equivalent to omit the latter
closure requirement on $N$. One may assume that the embedding $j$ comes from the ultrapower by an $M$-extender of length $|V_{\kappa+m-1}|$.

\Theorem. Every strongly unfoldable cardinal is totally indescribable and a limit of totally indescribable cardinals. Indeed, if $\kappa$ is merely
$(\kappa+\omega)$-strongly unfoldable, then the collection of totally indescribable cardinals below $\kappa$ is
stationary.\label{UnfoldablesAreIndescribable}

\Proof: Every strongly unfoldable cardinal is totally indescribable because the $(\kappa+m)$-strong unfoldability embeddings witness $\Pi^m_n$
indescribability for any $n$. Suppose that $\kappa$ is $(\kappa+\omega)$-strongly unfoldable and fix any club set $C\of\kappa$. Let $M$ be any
$\kappa$-model of set theory with $C\in M$. Since $M$ and $V$ agree up to $\kappa$, they agree on the totally indescribable cardinals below $\kappa$.
Fix any $\theta$-strong unfoldability embedding $j:M\to N$ for any $\theta\geq\kappa+\omega$. Thus, $V_{\kappa+\omega}\of N$. It follows that $N$ has
all $\kappa$-models and the extender embeddings on them to witness the $\Pi^m_n$ indescribability of $\kappa$, and consequently $\kappa$ is totally
indescribable in $N$. Since $\kappa$ is a limit point of $j(C)$, it follows that $\kappa\in j(C)$ and so $j(C)$ meets the set of totally
indescribable cardinals in $N$. By elementarity, $C$ contains some totally indescribable cardinals in $M$ and hence in $V$. So the set of such
cardinals is stationary.\QED

A cardinal $\kappa$ is {\df subtle} if for any closed unbounded set $C\of\kappa$ and any sequence $\<A_\alpha\st\alpha\in C>$ with
$A_\alpha\of\alpha$, there is $\alpha<\beta$ in $C$ with $A_\alpha=A_\beta\intersect\alpha$. It is not difficult to see that any such cardinal is
strongly inaccessible. We are grateful to Ralf Schindler for pointing out the following.

\Theorem. If $\kappa$ is subtle, then the set of cardinals $\gamma<\kappa$ that are $\ltkappa$-strongly unfoldable is stationary. Hence, $V_\kappa$
has a stationary proper class of strongly unfoldable cardinals.\label{SubtleIsStationaryLimit}

\Proof: Suppose not, so there is a closed unbounded set $C\of\kappa$ containing no $\ltkappa$-unfoldable cardinals. For each $\gamma$ in $C$, there
is some least $\lambda<\kappa$ such that $\gamma$ is not $\lambda$-strongly unfoldable, and by thinning $C$ we may assume that $\lambda$ is less than
the next element of $C$, and also that $\gamma$ is a Beth fixed point. Since $\gamma$ is not $\lambda$-strongly unfoldable, there is a transitive
structure $M_\gamma$ of size $\gamma$, with $\gamma\in M_\gamma$ and $V_\gamma\of M_\gamma$, having no embedding $j:M_\gamma\to N$ with
$\cp(j)=\gamma$, $j(\gamma)>\lambda$ and $V_\lambda\of N$. Since $M_\gamma$ has size $\gamma$, there is a relation $E_\gamma$ on $\gamma$ such that
$\<\gamma,E_\gamma>\iso \<M_\gamma,{\in}>$. The isomorphism $\pi_\gamma$ witnessing this is exactly the Mostowski collapse of $\<\gamma,E_\gamma>$.
We may assume that $\pi_\gamma(0)=\gamma$. Let $A_\gamma$ be a subset of $\gamma$ coding $E_\gamma$, the elementary diagram of $\<\gamma,E_\gamma>$
and the map $\pi_\gamma^{-1}\restrict\gamma$.

Since $\kappa$ is subtle, there is $\gamma<\delta$ in $C$ with $A_\gamma=A_\delta\intersect\gamma$. Define $j:M_\gamma\to M_\delta$ by
$j=\pi_\delta\compose\pi_\gamma^{-1}$. Observe that $j(\gamma)=\pi_\delta(\pi_\gamma^{-1}(\gamma))=\pi_\delta(0)=\delta$. Also, if $\alpha<\gamma$,
then because $A_\gamma$ and $A_\delta$ agree up to $\gamma$, it follows that $\pi_\gamma^{-1}\restrict\gamma=\pi_\delta^{-1}\restrict\gamma$, and so
$j(\alpha)=\alpha$. So the critical point of $j$ is $\gamma$. The map is elementary, since if $M_\gamma\satisfies\varphi[x]$ where
$x=\pi_\gamma(\alpha)$, then $\varphi(\alpha)$ is in the elementary diagram of $\<\gamma,E_\gamma>$, and so it is also in the elementary diagram of
$\<\delta,E_\delta>$, which means $M_\delta\satisfies\varphi[j(x)]$. Since $V_\delta\of M_\delta$, we have found a $\delta$-strong unfoldability
embedding for $M_\gamma$, contradicting our assumptions.\QED

Strongly unfoldable cardinals, like strong cardinals, have canonical extender embeddings, which have allowed for the borrowing of many techniques
from the strong cardinal context. Lemma \ref{Closure}, however, shows that they also exhibit a supercompactness-like nature, allowing us to borrow
techniques from the supercompact cardinal context as well. The fact that these miniature versions of strong and supercompact cardinals are equivalent
(see Corollary \ref{MiniatureStrongSupercompact}) is both interesting and surprising, so we include both characterizations.

\Lemma. If $\kappa$ is strongly $\theta$-unfoldable, then for any $\kappa$-model of set theory $M$ there is a strong $\theta$-unfoldability embedding
$j:M\to N$ such that every object in $N$ has the form $j(f)(\alpha)$ for some $\alpha<\beth_\theta$ and $f\in M\intersect M^\kappa$.\label{Canonical
Form Lemma}

\Proof: Suppose $j:M\to N$ is any strong $\theta$-unfoldability embedding. Let $X=\set{j(f)(\alpha)\st f\in M\cap M^\kappa \And
\alpha<\beth_\theta}$. By verifying the Tarski-Vaught criterion, one can easily see that $X\elesub N$, and in fact $X$ is the (smallest) Skolem hull
of $\ran(j)$ with the elements of $\beth_\theta$. If $j_0=\pi\compose j$, where $\pi:X\iso N_0$ is the Mostowski collapse of $X$, then $j_0:M\to N_0$
is an elementary embedding with critical point $\kappa$, and $j_0(\kappa)=\pi(j(\kappa))\geq\sup\pi\image\theta\geq\theta$. Let
$s=\<a_\alpha\st\alpha<\kappa>$ be an enumeration of $V_\kappa$ in $M$ with the property that the objects of any $V_\beta$ are enumerated before
$\beth_\beta$. It follows that every object in $V_\theta=(V_\theta)^N$ is $j(s)(\alpha)$ for some $\alpha<\beth_\theta$, and so $V_\theta\of X$. The
elements of $V_\theta$ are therefore fixed by the collapse $\pi$, and so $V_\theta\of N_0$. Finally, every object in $N_0$ has the form
$\pi(j(f)(\alpha))=j_0(f)(\pi(\alpha))=j_0(f)(\alpha)$. So $j_0:M\to N_0$ has the desired canonical form.\QED

These embeddings are obtained by the ultrapowers by $M$-extenders, in analogy with strong cardinals. One could just as easily insist alternatively
that every element of $N$ had the form $j(f)(a)$ for some $f\in M\intersect M^\kappa$ and $a\in V_\theta$. Next we turn to the supercompactness
nature of strongly unfoldable cardinals.

\Lemma. If $\kappa$ is $(\theta+1)$-strongly unfoldable, then for every $\kappa$-model of set theory $M$ there is a $\beth_\theta$-closed
$(\theta+1)$-strong unfoldability embedding $j:M\to N$. That is, $\cp(j)=\kappa$, $j(\kappa)>\theta$ and \begin{enumerate} \item $V_{\theta+1}\of N$,
\item $N^{\beth_\theta}\of N$ and \item $|N|=\beth_{\theta+1}$. \end{enumerate} \label{Closure}

\Proof: Suppose that $\kappa$ is $(\theta+1)$-strongly unfoldable and $M$ is a $\kappa$-model of set theory. Since $M$ has size $\kappa$, there is a
relation $E$ on $\kappa$ such that there is an isomorphism $\pi:\<\kappa,E>\iso\<M,{\in}>$. This isomorphism $\pi$ is necessarily the Mostowski
collapse of the structure $\<\kappa,E>$, and so it and the structure $M$ are uniquely determined by $E$. By collapsing the Skolem hull of $E$ in some
very large $V_\lambda$, we may construct a $\kappa$-model $\Mbar$ such that $E\in \Mbar$. It follows that $M\in \Mbar$ as well. By the
$(\theta+1)$-strong unfoldability of $\kappa$, there is an embedding $j:\Mbar\to\Nbar$ with critical point $\kappa$ such that $j(\kappa)>\theta$ and
$V_{\theta+1}\of\Nbar$. Since $E=j(E)\intersect\kappa\cross\kappa$, it follows that $E\in\Nbar$ and hence also $M\in\Nbar$. Let $j_0=j\restrict M$
and $N=j(M)$, so that $j_0:M\to N$. The critical point of $j_0$ is still $\kappa$ and $j_0(\kappa)=j(\kappa)>\theta$, so this is a
$(\theta+1)$-unfoldability embedding. Since $M^\ltkappa\of M$ and $\Mbar^\ltkappa\of\Mbar$, it follows that $M$ and $\Mbar$ agree up to $V_\kappa$.
Thus, $j(M)$ and $\Nbar$ agree up to $(V_{j(\kappa)})^\Nbar$. Since $\theta+1< j(\kappa)$, this includes the true $V_{\theta+1}$ and so $j_0:M\to N$
is a $(\theta+1)$-strong unfoldability embedding.

Now we use Hauser's trick from his treatment of indescribable embeddings in \cite{Hauser1991:IndescribableCardinals}. By elementarity, if $x\in M$ is
coded by some $\alpha<\kappa$ with respect to $E$, then $j(x)$ is coded by $j(\alpha)$ with respect to $j(E)$. That is, if $x=\pi(\alpha)$ for some
$\alpha<\kappa$, then $j(x)=j(\pi)(j(\alpha))$. But since $j_0(x)=j(x)$ and $j(\alpha)=\alpha$, this means that $j_0(x)=y$ if and only if $x$ is
coded by the same ordinal with respect to $E$ that $y$ is coded by with respect to $j(E)$. Consequently, the embedding $j_0$ is definable in $\Nbar$
from $E$ and $j(E)$ and their respective Mostowski collapses. Thus, the entire embedding $j_0:M\to N$ is an element of $\Nbar$.

Furthermore, since $\Mbar$ knows that $M^\ltkappa\of M$, it follows that $\Nbar$ satisfies $N^{<j(\kappa)}\of N$. In particular,
$N^{\beth_{\theta+1}}\of N$ in $\Nbar$. By a classical downward L\"owenheim-Skolem argument in $\Nbar$, we may build a set $X\elesub N$ such that
$\ran(j_0)\of X$, $X^{\beth_\theta}\of X$ and $|X|=\beth_{\theta+1}$ in $\Nbar$. We observe that any function from $\beth_\theta$ to $V_{\theta+1}$
can be canonically coded by an element of $V_{\theta+1}$, and so because $\Nbar$ has the true $V_{\theta+1}$, if $\Nbar$ thinks that a set of size
$\beth_{\theta+1}$, such as $X$, is closed under $\beth_\theta$ sequences, then $\Nbar$ is correct about this. So actually we know that
$X^{\beth_\theta}\of X$ in $V$. Let $\pi^*:X\iso N^*$ be the Mostowski collapse of $X$ and let $j^*=\pi^*\compose j_0$, so that $j^*:M\to N^*$. Since
$X^{\beth_\theta}\of X$, it follows that $V_\theta$ and all its subsets are in $X$. This implies $V_{\theta+1}\of X$. These objects will therefore be
fixed by $\pi^*$, and so $V_{\theta+1}\of N^*$ as well. Since $X$ has size $\beth_{\theta+1}$ and is isomorphic to $N^*$, it follows that $N^*$ also
has size $\beth_{\theta+1}$. Since $X^{\beth_\theta}\of X$ and $X\iso N^*$, it follows that $N^{\beth_\theta}\of N$. And finally, since
$\pi^*(\theta)=\theta$, it follows that $j^*(\kappa)=\pi^*(j_0(\kappa))\geq\pi^*(\theta+1)=\theta+1$. Thus, the embedding $j^*:M\to N^*$ is a
$\beth_\theta$-closed $(\theta+1)$-strong unfoldability embedding whose target is of size $\beth_{\theta+1}$, as we desired.\QED

In particular, if the \GCH\ holds at $\lambda=\beth_\theta$, then the lemma produces $\lambda$-strong unfoldability embeddings $j:M\to N$ such that
$N^\lambda\of N$ and $|N|=\lambda^\plus$. The same idea applies to the case when $\theta$ is a limit ordinal, but here we get only the $\cof(\theta)$
closure of $N$:

\Lemma. If $\kappa$ is $\theta$-strongly unfoldable for a limit ordinal $\theta$, then for any $\kappa$-model of set theory $M$ there is a
$\theta$-strong unfoldability embedding $j:M\to N$ with $\cp(j)=\kappa$, $j(\kappa)>\theta$ and \begin{enumerate} \item $V_\theta\of N$, \item
$N^{<\cof(\theta)}\of N$ and \item $|N|=\beth_\theta$. \end{enumerate} \label{LimitClosure}

\Proof: We proceed here just as above. Given $M$, we find $E$ and $\Mbar$ as above, and fix $j:\Mbar\to\Nbar$ a $\theta$-strong unfoldability
embedding. Once again, $j_0=j\restrict M$ is in $\Nbar$, and $\Nbar$ sees that $N=j(M)$ is closed under $<j(\kappa)$ sequences. By a
L\"owenheim-Skolem argument, we find $X\elesub N$ in $\Nbar$ such that $V_\theta\of X$, $X^{<\cof(\theta)}\of X$ in $\Nbar$ and $|X|=\beth_\theta$.
Since $\Nbar$ is correct about all the $<\cof(\theta)$ sequences in any set of size $\beth_\theta$, it follows that $X^{<\cof(\theta)}\of X$ in $V$
as well. And so if $N^*$ is the Mostowski collapse of $X$, it follows again that the corresponding embedding $j^*:M\to N^*$ is as we desired.\QED

\Corollary. The following are equivalent:
\begin{enumerate}
 \item The cardinal $\kappa$ is strongly unfoldable. That is, for every $\kappa$-model $M$ and ordinal $\theta$ there is an embedding
        $j:M\to N$ with $\cp(j)=\kappa$, $j(\kappa)>\theta$ and $V_\theta\of N$.
 \item For every $\kappa$-model $M$ and ordinal $\theta$ there is an embedding $j:M\to N$ with $\cp(j)=\kappa$,
        $j(\kappa)>\theta$ and $N^\theta\of N$.
\end{enumerate}\label{MiniatureStrongSupercompact}

So the miniature versions of strong and supercompact cardinals are equivalent. Most arguments in the literature involving strongly unfoldable
cardinals have until now used the strongness-like extender embeddings, but the supercompactness-like characterization allows for simplifications,
which we will make use of in our main argument.

\Section Warming-Up with Just One Sequence

In this section, we present the one-step forcing for destroying a single potential $\Diamond_\kappa(\REG)$ sequence while preserving the strong
unfoldability of $\kappa$.

\Definition. A {\df potential $\Diamond_\kappa$ sequence} is a sequence $\vec A=\<A_\alpha\st\alpha<\kappa>$ such that $A_\alpha\of\alpha$ for all
$\alpha<\kappa$. For any such sequence, the {\df diamond killing poset} $\Q_{\vec A}$ has conditions $(s,c)$ such that $c$ is a closed bounded subset
of $\kappa$, $s\of\sup(c)$ and $A_\alpha\neq s\intersect\alpha$ for every $\alpha\in c\intersect\INACC$. These conditions are ordered by
end-extension above $\sup(c)$, that is: $(s',c')\leq(s,c)$ if and only $s'\intersect\sup(c)=s$ and $c'\intersect(\sup(c)+1)=c$.

Potential $\Diamond_\kappa$ sequences are sometimes also referred to as `ladder systems', but the present terminology seems more evocative in our
case. The diamond killing poset clearly adds a set $S\of\kappa$ and  a club set $C\of\kappa$ such that $A_\alpha\neq S\intersect\alpha$ for any
$\alpha\in C\intersect\INACC$, meaning that it has killed off $\vec A$ as a witness for $\Diamond_\kappa(\REG)$.

\Lemma. Suppose $\kappa$ is inaccessible and $\vec A$ is a potential $\Diamond_\kappa$ sequence. Then: \begin{enumerate} \item The diamond killing
poset $\Q_{\vec A}$ preserves all cardinals and cofinalities. \item The diamond killing poset $\Q_{\vec A}$ adds no bounded sets to $\kappa$. \item
For every $\gamma<\kappa$, the diamond killing poset $\Q_{\vec A}$ has an open dense subset that is $\leqgamma$-directed closed. \item The diamond
killing poset $\Q_{\vec A}$ is $\kappa$-centered and thus satisfies the $\kappa^\plus$ chain condition. \item In the diamond killing extension
$V^{\Q_{\vec A}}$, or any further forcing extension, the sequence $\vec A$ does not witness $\Diamond_\kappa(\REG)$.
\end{enumerate}\label{KillingProperties}

\Proof: Statements 1 and 2 follow easily from 3 and 4. Statement 4 is immediate, because $\Q_{\vec A}$ has size $\kappa$. Statement 5 is also
immediate from our observation that $\Q_{\vec A}$ adds a set $S$ and a club set $C\of\kappa$ such that $A_\alpha\neq S\intersect\alpha$ for any
$\alpha\in C\intersect\INACC$. So it remains only to prove statement 3. For this, fix any $\gamma<\kappa$ and let $D$ be the set of conditions
$(s,c)$ such that $c\setminus\gamma$ is nonempty. This is certainly an open dense subset of $\Q_{\vec A}$, since we can strengthen any condition to
one mentioning ordinals above $\gamma$, and further strengthening only adds more elements to $c$. To see that $D$ is $\leqgamma$-directed closed,
suppose that $\set{(s_\alpha,c_\alpha)\st\alpha<\gamma}$ is a directed family in $D$. Let $s=\Union_\alpha s_\alpha$ and $c=\Union_\alpha c_\alpha$
be the limiting values of the coherent sets, and let $\bar c=c\union\singleton{\sup(c)}$ be the closure of $c$. Certainly
$s\intersect\delta=s_\alpha\intersect\delta$ for any $\delta\in c_\alpha\intersect\INACC$ for large enough $s_\alpha$, so $(s,\bar c)$ satisfies the
required nonreflection property $s\intersect\delta\not=A_\delta$ for any $\delta\in c\intersect\INACC$. The top point $\sup(c)$ itself cannot be
inaccessible, because it is above $\gamma$ but has cofinality at most $\gamma$. So $(s,\bar c)$ is a condition in $D$, and it is clearly below every
$(s_\alpha,c_\alpha)$. Statement 3 therefore holds.\QED

The condition $(s,\bar{c})$ constructed as above will be referred to as the limit of $\set{(s_\alpha,c_\alpha)\st\alpha<\gamma}$, and denoting
$q_\alpha=(s_\alpha,c_\alpha)$ and $q=(s,\bar{c})$, we may also say that $q=\inf_\alpha q_\alpha$.

While the diamond killing poset $\Q_{\vec A}$ kills off the sequence $\vec A$ as a witness for $\Diamond_\kappa(\REG)$, there are of course many
other sequences. In order to kill $\Diamond_\kappa(\REG)$ fully, one expects to iterate this forcing, anticipating and killing all the potential
diamond sequences that might arise. And with a wrinkle, this is how the proof of our main theorem will proceed.

Actually, if all one wants to do is kill the one sequence $\vec A$, there are much easier ways to do it. For example, by merely adding a single Cohen
real one preserves the strong unfoldability of $\kappa$ (and all other large cardinals) while destroying all diamond sequences in the ground model,
because they do not anticipate the new real. But clearly we cannot hope to iterate this forcing to achieve our purposes. Rather, the diamond killing
poset $\Q_{\vec A}$ is the natural choice for iterating when one wants not to add bounded subsets to $\kappa$.

Before we consider iterations of the diamond killing poset for different potential diamond sequences, we will analyze mere products of $\Q_{\vec A}$
with itself (where we kill the same sequence $\vec A$ on every coordinate; this is the wrinkle we mentioned above). Let $\Pi_{\kappa^\plus}\Q_{\vec
A}$ be the $\ltkappa$-support product of $\kappa^\plus$ many copies of $\Q_{\vec A}$. It may seem odd to consider forcing that kills the same
sequence $\vec A$ so many times, since one might expect that killing it once would be enough. But in order to preserve the strong unfoldability of
$\kappa$ in the extension, this many-times-over death of $\vec A$ as a witness for $\Diamond_\kappa(\REG)$ is what the argument seems to require. Our
argument, specifically, will make critical use of the fact that we have $\kappa^\plus$ many different generic filters for $\Q_{\vec A}$ available to
us. Hauser's argument \cite{Hauser92:IndescribablesWithoutDiamond} for the $\Pi^m_n$ indescribable cardinals exhibits the same repetitive killing
feature, where he makes similar crucial use of the fact that every potential diamond sequence is killed unboundedly often during his iteration. One
difference in presentation is that while Hauser blends all his repetitions together in one long iteration, we have isolated the copies of $\Q_{\vec
A}$ all in one product. This will allow us to illustrate the general technique as a warm-up in Theorem \ref{WarmUp}, in the simplified case of just
one potential diamond sequence $\vec A$.

\Lemma. Suppose $\kappa$ is inaccessible and $\vec A$ is a potential diamond sequence. Let $\Pi_{\kappa^\plus}\Q_{\vec A}$ be the $\ltkappa$-support
product of the diamond killing poset $\Q_{\vec A}$. Then: \begin{enumerate} \item The forcing $\Pi_{\kappa^\plus}\Q_{\vec A}$ preserves all cardinals
and cofinalities. \item The forcing $\Pi_{\kappa^\plus}\Q_{\vec A}$ adds no bounded sets to $\kappa$. \item For every $\gamma<\kappa$, the forcing
$\Pi_{\kappa^\plus}\Q_{\vec A}$ has a dense subset that is $\leqgamma$-directed closed. \item The forcing $\Pi_{\kappa^\plus}\Q_{\vec A}$ is
$\kappa$-centered and thus satisfies the $\kappa^\plus$ chain condition. \item In the extension $V^{\Pi_{\kappa^\plus}\Q_{\vec A}}$, or any further
forcing extension, the sequence $\vec A$ does not witness $\Diamond_\kappa(\REG)$.
\end{enumerate}\label{ProductKillingProperties}

\Proof: Again, statements 1 and 2 follow from 3 and 4. Statement 5 is immediate from Lemma \ref{KillingProperties}, since
$V^{\Pi_{\kappa^\plus}\Q_{\vec A}}$ is an extension of $V^{\Q_{\vec A}}$. For statement 3, we consider the set $D$ consisting of conditions that are
in the $\leqgamma$-closed dense subset of $\Q_{\vec A}$ on each coordinate in their support.  That is, $D$ consists of conditions in the product such
that every coordinate is either trivial or else mentions ordinals at least as large as $\gamma$. This is easily seen to be dense and
$\leqgamma$-directed closed by the same argument as in Lemma \ref{KillingProperties}, so 3 holds. Statement 4 follows from the fact that
$\kappa^\ltkappa=\kappa$ and the poset is the $\ltkappa$-support product of $\kappa$-centered forcing.\QED

We warn the reader that in Lemma \ref{ProductKillingProperties}, the $\leqgamma$-directed closed dense set we produce in statement 3 is not open, as
in Lemma \ref{KillingProperties}, because one can always extend a condition by adding an element to the support and placing a condition there that
does not jump up to $\gamma$.

For convenience, our main construction employs the lottery preparation of \cite{Hamkins2000:LotteryPreparation}. While one could avoid this by
tailoring the iteration below $\kappa$ to ensure the desired reflectivity, the point of the lottery preparation is that such efforts are unnecessary:
one simply works below a condition opting for the desired forcing at stage $\kappa$, and the generic filter is thereby forced to have the desired
reflectivity.

Let us quickly review the lottery preparation here. By forcing if necessary, we may assume that there is a function $f\from\kappa\to\kappa$ with the
strong unfoldability {\df Menas} property: for every $\kappa$-model $M$ and every ordinal $\theta$ there is a $\theta$-strong unfoldability embedding
$j:M\to N$ such that $j(f)(\kappa)>\theta$. Arguments in \cite{Hamkins2001:UnfoldableCardinals} show that such a function can be added by Woodin's
fast function forcing, while preserving the strong unfoldability of $\kappa$. And it is easy to see that the collapse of $X$ in Lemmas \ref{Closure}
and \ref{LimitClosure} will preserve the fact that $j(f)(\kappa)>\theta$, so we may take $j:M\to N$ to have the supercompactness-like form of those
lemmas. The {\df lottery preparation} relative to $f$ is the Easton support $\kappa$ iteration $\P$ which at every stage $\gamma<\kappa$, provided
that $\gamma\in\dom(f)$ and $f\image\gamma\of\gamma$, forces with the lottery sum of all $\Q\in H(f(\gamma))$ having for every $\beta<\gamma$ a dense
subset that is $\leqbeta$-strategically closed. (The {\df lottery sum} $\oplus{\cal A}$ of a collection of posets $\cal A$, also commonly called
side-by-side forcing, is the poset $\{\,\<\Q,p>\st p\in\Q\in{\cal A}\,\}\union\singleton{\one}$, ordered with $\one$ above everything and
$\<\Q,p>\leq\<\Q',p'>$ when $\Q=\Q'$ and $p\leq_\Q p'$. The generic filter in effect selects a `winning' poset from $\cal A$ and then forces with
it.) By further restricting the lotteries to include only forcing notions that preserve the \GCH, one obtains the {\df \GCH-preserving} lottery
preparation, and this preserves the \GCH. The thrust of \cite{Hamkins2000:LotteryPreparation} is that if $j$ is an embedding with critical point
$\kappa$, then the lottery sum at stage $\kappa$ in $j(\P)$ includes all the desired posets, and so by working below a condition opting for the
correct poset in that lottery, one avoids the need for a Laver function. Arguments in \cite{Hamkins2001:UnfoldableCardinals} show that the lottery
preparation of a strongly unfoldable cardinal $\kappa$ preserves the strong unfoldability of $\kappa$, and makes the strong unfoldability of $\kappa$
indestructible, for example, by the forcing $\add(\kappa,1)$, among others. And these arguments are easily adapted to the \GCH-preserving lottery
preparation.

We turn now to the Warm-Up Theorem, which shows how to kill off a single potential diamond sequence $\vec A$, while preserving the strong
unfoldability of $\kappa$ and adding no new bounded sets to $\kappa$. Later, in the Main Theorem, we will simply iterate this forcing so as to
anticipate all possible sequences $\vec A$, thereby forcing the negation of $\Diamond_\kappa(\REG)$ while preserving the strong unfoldability of
$\kappa$.

\break

\Theorem The Warm-Up Theorem. Suppose that $\kappa$ is strongly unfoldable in $V$, the \GCH\ holds and $V[G]$ is the \GCH-preserving lottery
preparation of $\kappa$ relative to $f\from\kappa\to\kappa$. Then for any potential diamond sequence $\vec A$ in $V[G]$, the cardinal $\kappa$
remains strongly unfoldable after further forcing with $\Pi_{\kappa^\plus}\Q_{\vec A}$.\label{WarmUp}

\Proof: Suppose that $g\of\Pi_{\kappa^\plus}\Q_{\vec A}$ is $V[G]$-generic for the $\ltkappa$-support $\kappa^\plus$ product of the diamond killing
poset $\Q_{\vec A}$ corresponding to $\vec A=\<A_\gamma\st\gamma<\kappa>$. This forcing preserves the \GCH\ because it adds no bounded subsets to
$\kappa$ and has size $\kappa^\plus$. The generic filter $g$ is determined by the generic objects $(S_\xi,C_\xi)$ added by the copy of $\Q_{\vec A}$
at each coordinate $\xi<\kappa^\plus$. Thus, $C_\xi\of\kappa$ is a club set, $S_\xi\of\kappa$ and for every $\gamma\in C_\xi\intersect\INACC$ we have
$S_\xi\intersect\gamma\neq A_\gamma$.

We claim that $\kappa$ is strongly unfoldable in $V[G][g]$. If not, then there is a condition $p\in g$, an ordinal $\theta>\kappa$ and a name $\Bdot$
for a subset of $\kappa$ such that in $V[G]$, the condition $p$ forces that $\Bdot$ cannot be placed into a $\kappa$-model having a $\theta$-strong
unfoldability embedding. By increasing $\theta$ if necessary, we may assume that $\theta=\beth_\theta$. Because the forcing is $\kappa^\plus$-c.c.,
we may assume that the name $\Bdot$ has hereditary size at most $\kappa$ and depends on at most $\kappa$ many coordinates of the product forcing.

Let $M$ be a $\kappa$-model of set theory in $V$, large enough to contain the things (of size $\kappa$) on which we are focused: $\kappa$, $f$, $\P$
and names for $\vec A$, $\Bdot$ and $p$, so that these latter items are elements of $M[G]$. Because $\kappa$ is strongly unfoldable in $V$ and the
function $f$ has the Menas property, by Lemma \ref{Closure} there is a $(\theta+1)$-strong unfoldability embedding $j:M\to N$ in $V$ with critical
point $\kappa$ such that $j(f)(\kappa)>\theta$, $N^\theta\of N$, $V_{\theta+1}\of N$ and $|N|=\theta^\plus$.  Since $V_\theta\of N$, the model $N$
has the same $\kappa^\plus$ as $V$, and so the full forcing $\Pi_{\kappa^\plus}\Q_{\vec A}$ appears as one of the choices in the stage $\kappa$
lottery of $j(\P)$. Below the condition opting for this poset, the forcing $j(\P)$ factors as $\P*(\Pi_{\kappa^\plus}\Q_{\vec A})*\Ptail$, where
$\Ptail$ is the forcing from stages $\kappa+1$ up to $j(\kappa)$. Note that $G*g$ is $N$-generic for the first $\kappa+1$ many stages of this
forcing. Since $j(f)(\kappa)>\theta$, the next nontrivial stage of forcing in $\Ptail$ is beyond $\theta$, and so $\Ptail$ has a
$\leqtheta$-strategically closed dense set in $N[G][g]$. Furthermore, there are only $\theta^\plus$ many elements in $N$ altogether, so in $V[G][g]$
we may list the dense subsets of $\Ptail$ in $N[G][g]$ in a $\theta^\plus$-sequence. Since $N^\theta\of N$ in $V$ and the forcing $G*g$ is
$\kappa^\plus$-c.c., it follows that $N[G][g]^\theta\of N[G][g]$ in $V[G][g]$. We may therefore construct by diagonalization in $V[G][g]$ a
descending sequence of conditions in $\Ptail$ that meet every dense set in $N[G][g]$. These generate therefore an $N[G][g]$-generic filter
$\Gtail\of\Ptail$, and we may lift the embedding to $j:M[G]\to N[j(G)]$, where $j(G)=G*g*\Gtail$. Since $\vec A$ has a name in $M$, it follows that
$\vec A\in M[G]$ and so $j(\vec A)$ is a $j(\kappa)$-sequence in $N[j(G)]$.

Consider now the forcing $(\Pi_{\kappa^\plus}\Q_{\vec A})^{M[G]}$, which is the same as $\Pi_\delta\Q_{\vec A}$ in $V[G]$, where
$\delta=(\kappa^\plus)^M$. It follows that $g\restrict\delta$ is $M[G]$-generic for this forcing. The lucky case occurs if it happens that
$S_\alpha\neq j(\vec A)(\kappa)$ for all $\alpha<\delta$. In this case, we use $j\image\delta\in N$ to build the (master) condition $q\in
j(\Pi_\delta\Q_{\vec A})$ having support $j\image\delta$, with the pair $(S_\alpha,C_\alpha\union\singleton{\kappa})$ at coordinate $j(\alpha)$.
First of all, this object is an element of $N[j(G)]$ because $g\in N[j(G)]$ and $j\image\delta\in N$. It is a condition in $j(\Pi_\delta\Q_{\vec A})$
because it has the correct support, it satisfies all the required non-reflection conditions below $\kappa$, and it satisfies the non-reflection
condition at stage $\kappa$ by the lucky property that $j(S_\alpha)\intersect\kappa=S_\alpha\neq j(\vec A)(\kappa)$. This condition is called a
master condition because it is stronger than every element of $j\image (g\restrict\delta)$. Below this (or any) condition, the forcing
$j(\Pi_\delta\Q_{\vec A})$ has a $\leqtheta$-strategically closed dense set, so using the fact that $N[j(G)]^\theta\of N[j(G)]$ in $V[G][g]$ we may
once again diagonalize to build a $N[j(G)]$-generic filter containing $q$. The embedding therefore lifts in $V[G][g]$ to $j:M[G][g\restrict\delta]\to
N[j(G)][j(g\restrict\delta)]$. And since $V_{\theta+1}\of N$ and $G*g\in N[j(G)]$, it follows that $(V[G][g])_{\theta+1}\of N[j(G)]$, making the lift
a strong $(\theta+1)$-unfoldability embedding in $V[G][g]$. Finally, since $\Bdot\in M[G]$, the value of $\Bdot_g$ depends only on the first $\delta$
many coordinates of $g$, so $B=\Bdot_g=\Bdot_{g\restrict\delta}\in M[G][g\restrict\delta]$, which contradicts our assumption that there was no
$\theta$-strongly unfoldable embedding for $B$.

In general, however, we may not be in the lucky case (see Remark \ref{RemarkOnLuckyCase}), and our task is a bit harder. Suppose the unlucky case
occurs, and $S_{\alpha_0}=j(\vec A)(\kappa)$ for some $\alpha_0<\delta$. Since the $S_\alpha$ are distinct, this occurs for at most one
$\alpha_0<\kappa^\plus$. Consider the original condition $p$, forcing that $\Bdot$ named a counterexample, and it's value $p(\alpha_0)$ on coordinate
$\alpha_0$. It is dense that there is some other coordinate $\alpha_1\geq\delta$ such that the generic object on coordinate $\alpha_1$ extends
$p(\alpha_0)$ (after all, any given condition in $\Q_{\vec A}$ will occur in the product generic on unboundedly many coordinates). Let $\pi^*$ be the
automorphism of $\Pi_{\kappa^\plus}\Q_{\vec A}$ determined by swapping coordinates $\alpha_0$ and $\alpha_1$. If $g^*=\pi^*\image g$, then our choice
of $\alpha_1$ ensures that $p\in g^*$. Furthermore, because the automorphism is in the ground model, it follows that $V[G][g]=V[G][g^*]$. And
finally, since $S_{\alpha_1}\neq S_{\alpha_0}=j(\vec A)(\kappa)$, the generic $g^*$ is in the lucky case when it comes to lifting the embedding
$j:M[G]\to N[j(G)]$ that we have already constructed so far. So as above we may construct a master condition below $j\image g^*\restrict\delta$ and
an $M[j(G)]$-generic filter containing it, thereby lifting the embedding to $j:M[G][g^*\restrict\delta]\to N[j(G)][j(g^*\restrict\delta)]$. Since we
arranged that $p\in g^*$ and $B^*=\Bdot_{g^*}\in M[G][g^*\restrict\delta]$, this contradicts once again our assumption that $p$ forced that $\Bdot$
named a set that could not be placed into a $\theta$-strong unfoldability embedding.\QED

\Remark. \label{RemarkOnLuckyCase} We emphasize that in the argument above, we are not claiming that every generic $g$ has an isomorphic copy $g^*$
placing $V[G][g^*]$ into the lucky case. This is actually false, because the sequence $\vec A$ could be defined so that $A_\gamma$ picks out the
first generic set $S_0\of\gamma$ added by the forcing $G$ at stage $\gamma$, if the forcing $G$ opted at stage $\gamma$ for a $\gamma^\plus$ product
of diamond killing posets. In this case, if one runs the construction from scratch in $V[G][g^*]$, we will have $j(G)=G*g^**\Gtail$ and consequently
$j(\vec A)(\kappa)=S_0^*$, the first set added by $g^*$, placing us in the unlucky case regardless of which $g^*$ we use. The point of the argument,
rather, is to use the original generic $g$ to build the partial lift $j:M[G]\to N[j(G)]$ with $j(G)=G*g*\Gtail$, thereby fixing the value of $j(\vec
A)(\kappa)$, and then to shift midstream to an isomorphic copy $g^*$ that puts us into the lucky case for the rest of the argument. The remarks at
the end of this article show that in general there can be no way for a pure lifting argument to succeed when one is trying to destroy
$\Diamond_\kappa(\REG)$, because such arguments will preserve the existence of a strong unfoldability Laver function, which implies
$\Diamond_\kappa(\REG)$.

\Section The Main Theorem

We will now iterate the one-step forcing of the previous section in order to kill off all the potential $\Diamond_\kappa(\REG)$ sequences that may
arise.

\Theorem Main Theorem. If $\kappa$ is strongly unfoldable, then there is a forcing extension preserving this in which $\Diamond_\kappa(\REG)$ fails
and the \GCH\ holds.

\Proof: The main result of \cite{Hamkins2001:UnfoldableCardinals} shows in part that by forcing if necessary, we may assume that the \GCH\ holds and
there is a function $f$ with the Menas property for strong unfoldability. Suppose $G\of\P$ is $V$-generic for the \GCH-preserving lottery preparation
of $\kappa$ with respect to $f$. Let $\Q$ be the $\ltkappa$-support forcing iteration which at each stage $\nu<\kappa^\plus$ forces with the
$\kappa^\plus$ product $\Pi_{\kappa^\plus}\Q_{\vec A^\nu}$ for some potential diamond sequence $\vec A^\nu=\<A^\nu_\gamma\st\gamma<\kappa>$, chosen
by a book-keeping function in such a way that all potential diamond sequences added by $\Q$ are handled by some stage. Such a book-keeping function
exists because the iteration $\Q$, being the $\ltkappa$-support iteration of $\kappa$-centered forcing at each stage, is $\kappa$-centered and hence
$\kappa^\plus$-c.c.; so any potential diamond sequence will appear by some stage of the forcing. And since any such stage of the forcing adds at most
$2^\kappa=\kappa^\plus$ such sequences, we can anticipate all the sequences that arise in one $\kappa^\plus$-iteration. Suppose $H\of\Q$ is
$V[G]$-generic for this forcing. The generic object $H$ is determined, therefore, by the sequence of sequences
$\<\<(S^\nu_\alpha,C^\nu_\alpha)\st\alpha<\kappa^\plus>\st\nu<\kappa^\plus>$, where $(S^\nu_\alpha,C^\nu_\alpha)$ are the sets added at the
$\alpha^\th$ coordinate of the product in the stage $\nu$ forcing. By the properties identified in Lemma \ref{ProductKillingProperties} for each step
of the iteration, it follows that the iteration altogether preserves all cardinals and cofinalities, adds no bounded subsets to $\kappa$, has for any
$\gamma<\kappa$ a $\leqgamma$-directed closed dense set and satisfies the $\kappa^\plus$-chain condition. Because the book-keeping function $\<\,\vec
A^\nu\st\nu<\kappa^\plus>$ eventually considers every potential diamond sequence in $V[G][H]$, killing $\vec A^\nu$ at stage $\nu$, it follows that
$V[G][H]\satisfies\neg\Diamond_\kappa(\REG)$. Finally, the \GCH-preserving lottery preparation preserves the \GCH\ to $V[G]$, and the remaining
forcing has size $\kappa^\plus$ and adds no bounded subsets to $\kappa$, so $V[G][H]$ satisfies the \GCH.

It remains to show that $\kappa$ remains strongly unfoldable in $V[G][H]$. Suppose towards contradiction that it does not. Then there is a condition
$p\in H$, an ordinal $\theta>\kappa$ and a name $\Bdot$ for a subset of $\kappa$ such that in $V[G]$ the condition $p$ forces that $\Bdot$ is not
captured by any $\theta$-strong unfoldability embedding. We may assume as in the warm-up argument that $\theta=\beth_\theta$ and that $\Bdot$ is a
nice name. By the chain condition and the support in the iteration $\Q$, it follows that $\Bdot$ mentions at most $\kappa$ many stages of $\Q$ and
within any such stage, mentions at most $\kappa$ many coordinates (recall that each stage $\nu$ of the forcing is the $\ltkappa$-support product
$\Pi_{\kappa^\plus}\Q_{\vec A^\nu }$ of the diamond killing poset $\Q_{\vec A^\nu}$). The name $\Bdot$ therefore has hereditary size at most
$\kappa$. Let $\lambda\muchgt\theta$ be an enormous cardinal such that $V_\lambda$ is a model of set theory, and construct $X\elesub V_\lambda$ of
size $\kappa$ in $V$ such that $X^\ltkappa\of X$ and $X$ contains all of the objects on which we are focused: $\kappa$, $\P$, $\Q$, $\Bdot$, $p$ and
the book-keeping function $\<\vec A^\nu \st\nu<\kappa^\plus>$. Let $M$ be the Mostowski collapse of $X$. Because $X^\ltkappa\of X$ and $X\elesub
V_\lambda$, it follows that $M$ is a $\kappa$-model of set theory. Furthermore, since the objects of hereditary size $\kappa$ are not affected by the
collapse of $X$ to $M$, we know that $\kappa$, $\P$, $p$ and $\Bdot$ are in $M$. Since $G$ is $V$-generic, it is also $M$-generic and we may form the
extension $M[G]$, the \GCH-preserving lottery preparation of $\kappa$ with respect to $f$. Since the book-keeping function is collapsed to the
initial segment of itself of length $(\kappa^\plus)^M$, the iteration $\Q^{M[G]}$ is the $(\kappa^\plus)^M$ iteration which at stage
$\nu<(\kappa^\plus)^M$ forces with the product of $(\kappa^\plus)^M$ many copies of the diamond killing poset $\Q_{\vec A^\nu }$. Thus, $\Q^{M[G]}$
is completely embedded in $\Q$, and the restriction of $H$ to the appropriate domains will be $V[G]$-generic for $\Q^{M[G]}$.

Let $j:M\to N$ be a $(\theta+1)$-strong unfoldability embedding in $V$ witnessing the Menas property of $f$, so that $j(f)(\kappa)>\theta$,
$N^\theta\of N$ in $V$, $V_{\theta+1}\of N$ and $|N|=\theta^\plus$. As in the warm up argument of Theorem \ref{WarmUp}, we may lift the embedding in
$V[G][H]$ to $j:M[G]\to N[j(G)]$, where $j(G)=G*H*\Gtail$. (This proceeded, as you recall, by opting for the full iteration $\Q$ in the stage
$\kappa$ lottery of $j(\P)$, using the generic object $H$ at this stage and building the remainder of the generic $\Gtail$ by diagonalization.) We
will now lift the embedding through the forcing $\Q^{M[G]}$. As in the warm-up argument, while the restriction of $H$ does provide an $M[G]$-generic
filter for this forcing, it may not admit a master condition that allows us to lift the embedding. As before, we will instead build a master
condition for a suitable isomorphic copy of $H$ that does admit a master condition.

\SubLemma. There is an automorphism $\pi:\Q\iso\Q$ and a condition $q\in j(\Q^{M[G]})$ such that if $H^*=\pi\image H$, then $p\in H^*$ and $q$ is
below every element of $j\image H^*\restrict\Q^{M[G]}$.

\Proof: The automorphism $\pi$ will be essentially an iteration of the kinds of automorphisms that arose in Theorem \ref{WarmUp}, which merely
swapped two coordinates in a large product. Here, at each stage $\nu<(\kappa^\plus)^M$, the automorphism $\pi$ will either do nothing to the stage
$\nu$ forcing or else swap two coordinates in the product at stage $\nu$. At stages beyond $(\kappa^\plus)^M$, the automorphism will do nothing extra
(though of course any automorphism of the earlier stages naturally induces a map in the names forming the later stages, which we intend to be carried
out without further remarks). The master condition $q$ will be the limit of a descending sequence of conditions $q_\beta$ for
$\beta<(\kappa^\plus)^M$ that we construct simultaneously along with $\pi$. Each of those proto-master conditions $q_\beta$ will be chosen from the
fixed $\leqkappa$-closed dense set $D\of j(\Q^{M[G]})$ consisting as in Lemma \ref{ProductKillingProperties} of those conditions which at any stage
and on  any coordinate in the product at that stage, are (with Boolean value 1) either trivial or else mention ordinals at least as large as
$\kappa$. Thus, $D$ is the iteration of the sets shown in Lemma \ref{ProductKillingProperties} to be $\leqkappa$-closed at each stage, and so it is
$\leqkappa$-closed in $N[j(G)]$, and the restriction of $D$ to any stage of forcing is also $\leqkappa$-closed.

Suppose now that the conditions $q_\beta$ have been chosen for all $\beta<\nu$, where $\nu<(\kappa^\plus)^M$, and that the action of the automorphism
$\pi$ has been specified on $\Q\restrict\nu$, the iteration of all the stages of forcing up to $\nu$. We assume inductively that $q_\beta\in
j(\Q^{M[G]})\restrict j(\beta)$, so that $q^*=\inf_{\beta<\nu}q_\beta\in j(\Q^{M[G]})\restrict\sup j\image\nu$. The lucky case occurs, as in Theorem
\ref{WarmUp}, when there is a condition $q^{**}$ below $q^*$ in $D\intersect j(\Q^{M[G]})\restrict j(\nu)$ forcing that $j(\vec A^\nu)_\kappa$ is not
equal to $S^\nu_\alpha$ for any $\alpha<(\kappa^\plus)^M$. In this case, we let $q_\nu=q^{**}\concat r$ where $r$ is the name of the condition at
stage $j(\nu)$ having support $j\image(\kappa^\plus)^M$ in the $j((\kappa^\plus)^M)$ product at that stage, and placing $(S^\nu_\alpha,\bar
C^\nu_\alpha)$ at coordinate $j(\alpha)$, where $\bar C^\nu_\alpha=C^\nu_\alpha\union\singleton{\kappa}$. That is,
$q_\nu(j(\nu),j(\alpha))=r(j(\alpha))=(S^\nu_\alpha,\bar C^\nu_\alpha)$. Note that this is a condition, because the support is smaller than
$j(\kappa)$, and the lucky condition guarantees exactly the nonreflection property we need at $\kappa\in\bar C^\nu_\alpha$. Also, since
$\kappa\in\bar C^\nu_\alpha$, this condition remains in the dense set $D$. Finally, in this lucky case, the automorphism does nothing to the stage
$\nu$ forcing.

In the unlucky case, there will be a condition $q^{**}\leq q^*$ in $D\intersect j(\Q^{M[G]})\restrict j(\nu)$ and an ordinal
$\alpha_0<(\kappa^\plus)^M$ such that $q^{**}$ forces $j(\vec A^\nu)_\kappa=S^\nu_{\alpha_0}$. As in Theorem \ref{WarmUp}, we may find an ordinal
$\alpha_1\geq(\kappa^\plus)^M$ and a still stronger condition $q^{***}$ forcing that the value of $p(\nu,\alpha_0)$ is included in
$(S^\nu_{\alpha_1},C^\nu_{\alpha_1})$. In this case, we let $q_\nu=q^{***}\concat r$ where again $r(j(\nu),j(\alpha))=(S^\nu_{\alpha},\bar
C^\nu_{\alpha})$ for $\alpha\neq\alpha_0$, but $r(j(\nu),j(\alpha_0))=(S^\nu_{\alpha_1},\bar C^\nu_{\alpha_1})$. Let the automorphism $\pi$ act on
the stage $\nu$ forcing by swapping the coordinates $\alpha_0$ and $\alpha_1$ in the product at stage $\nu$. This completes the recursive definition
of the $q_\beta$ and of $\pi$. Let $q$ be the limit of the conditions $q_\beta$ for $\beta<(\kappa^\plus)^M$. This limit exists because
$(\kappa^\plus)^M<\kappa^\plus$ and $D$ is $\leqkappa$-closed.

Let $H^*$ be the image of $H$ under the automorphism $\pi$. Observe that by the choice of $\alpha_1$ at stage $\nu$, it follows that
$H(\nu,\alpha_1)$ extends $p(\nu,\alpha_0)$, and so $H^*(\nu,\alpha_0)$ extends $p(\nu,\alpha_0)$. Since on the other coordinates, $H^*$ and $H$
agree, it follows that $p\in H^*$. Finally, since by construction the condition at $q(j(\nu),j(\alpha))$ for $\nu,\alpha<(\kappa^\plus)^M$ extends
every condition in the generic at $H^*(\nu,\alpha)$, it follows that $q$ is below every condition in $j\image H^*\restrict\Q^{M[G]}$, as desired.\QED

To complete the proof of the theorem, given the lemma, we proceed as in the warm-up of Theorem \ref{WarmUp}. Since $q$ is a master condition for the
forcing $H^*\restrict\Q^{M[G]}$, we can use the diagonalization technique to lift the embedding fully to $j:M[G][H^*\restrict\Q^{M[G]}]\to
N[j(G)][j(H^*\restrict\Q^{M[G]})]$. As before, this is a $(\theta+1)$-strong unfoldability embedding in $V[G][H]=V[G][H^*]$. If
$B=\Bdot_{H^*}=\Bdot_{H^*\restrict\Q^{M[G]}}$, then as before it follows that $B\in M[G][H^*\restrict\Q^{M[G]}]$, and so this embedding has captured
$B$. But the condition $p$ forced that $\Bdot$ named a set that could not be captured by any $\theta$-strong unfoldability embedding. This
contradicts the fact that $p\in H^*$.\QED

\Corollary. The existence of an unfoldable cardinal is equiconsistent with the existence of a strongly unfoldable cardinal $\kappa$ satisfying
$\neg\Diamond_\kappa(\REG)$, plus the \GCH.

The proof of the main theorem can be refined to give the following more local result.

\Theorem. If $\kappa$ is $(\theta+1)$-strongly unfoldable, then there is a forcing extension preserving this in which $\neg\Diamond_\kappa(\REG)$
holds.

\Proof: The point is that in the proof above, we do not really need to assume that $\theta=\beth_\theta$. The argument works perfectly well just
knowing that the \GCH\ holds at $\beth_\theta$. And one can force this by simply adding a Cohen subset to $\beth_{\theta}^\plus$. This forcing will
ensure that $\beth_{\theta+1}=\beth_\theta^\plus$, and since it is $\leq\beth_\theta$-closed, it will add no new $\kappa$-models. One can now run the
argument above, lifting the embedding $j:M\to N$, but modified to place the Cohen subset of $\beth_\theta^\plus$ in the stage $\kappa$ forcing of
$j(\P)$, as well as the $\kappa^\plus$ iteration $\Q$. Since $N$ has size $\beth_\theta^\plus$ and is $\leq\beth_\theta$-closed, the diagonalization
arguments are still able to produce the tail generic filters, so that $j$ lifts to $j:M[G]\to N[j(G)]$ with $j(G)=G*(H*A)*\Gtail$, where $H\of\Q$ is
the generic for the diamond killing iteration, $A$ is the Cohen subset of $\beth_\theta^\plus$ and $\Gtail$ is constructed by diagonalization. The
master condition and automorphism argument go through unchanged to produce $j:M[G][H]\to N[j(G)][j(H)]$. This will be a $(\theta+1)$-strong
unfoldability embedding because $V_{\theta+1}\of N$ and we placed $A*H$ at stage $\kappa$ in $j(\P)$. Since this embedding captures the set $B$,
whose name we placed into $M$ as above, we conclude that $\kappa$ is $(\theta+1)$-strongly unfoldable in $V[G][H][A]$. (Note that there is no need to
lift the embedding through the forcing to add $A$, because every subset of $\kappa$ is contained in a $\kappa$-model in $V[G][H]$.)\QED

\Section Further Observations

We will now show how to apply our analysis to separate two intimately connected combinatorial principles: the existence of an ordinal anticipating
Laver function versus a full set-anticipating Laver function, in the large cardinal contexts from weakly compact to strong unfoldable. The fact that
these are not equivalent answers an open question in \cite{Hamkins:LaverDiamond}.

\Definition.(\cite{Hamkins:LaverDiamond}) The {\df Laver Diamond principle} $\LDwc_\kappa$ for weak compactness is the assertion that there is a
function $\ell\from\kappa\to V_\kappa$ such that for any $A\in H(\kappa^\plus)$ and any $\kappa$-model $M$ containing $\ell$ and $A$ there is an
embedding $j:M\to N$ with critical point $\kappa$ such that $j(\ell)(\kappa)=A$. (An equivalent principle is investigated independently in
\cite{ShelahVaananan726:ExtensionsOfInfinitaryLogic} and \cite{Hellston2004:DiamondsOnLC}. The principle implies $\Diamond_\kappa(\REG)$ because any
club set $C\of\kappa$ can be placed along with $A$ into such an $M$, and so from $j(\ell)(\kappa)=A$ and $\kappa\in j(C)$ it follows that $\ell$
anticipates $A$ in $C\intersect\REG$.) The function $\ell$ is a mere {\df ordinal-anticipating} Laver function if for any $\kappa$-model $M$
containing $\ell$ and any ordinal $\alpha<\kappa^\plus$ there is an embedding $j:M\to N$ such that $j(\ell)(\kappa)=\alpha$. The Laver Diamond
principle $\LDsunf_\kappa$ for strong unfoldability asserts that there is a function $\ell\from\kappa\to V_\kappa$ such that for any ordinal
$\theta$, any $A\in V_\theta$ and any $\kappa$-model $M$ containing $\ell$, there is a $\theta$-strong unfoldability embedding $j:M\to N$ such that
$j(\ell)(\kappa)=A$. An {\df ordinal-anticipating} strong unfoldability Laver function is only required to anticipate every ordinal $\alpha<\theta$
with $j(\ell)(\kappa)=\alpha$. One can similarly define the Laver Diamond principle $\LDind_\kappa$ for indescribability and $\LDunf_\kappa$ for
unfoldability, as well as for many other large cardinals, along with the ordinal-anticipating versions of these principles.

\Theorem. If $V$ has a Laver Diamond $\LDsunf_\kappa$ function, then the extension $V[G][H]$ of the Main Theorem nearly does. Specifically, there is
a function $\ell^*\from\kappa\to V[G]_\kappa$ in $V[G][H]$ with the following properties for any $\theta$:
\begin{enumerate}
 \item For any $B\in P(\kappa)^{V[G]}$ and any $A\in V_\theta[G][H]$, there is a $\kappa$-model $M$ with $B,\ell^*\in M$ and a $\theta$-strong unfoldability
 embedding $j:M\to N$ with $j(\ell^*)(\kappa)=A$.
 \item For any $B\in P(\kappa)^{V[G][H]}$ and any $A\in V_\theta[G]$, there is a $\kappa$-model $M$ with $B,\ell^*\in M$ and a $\theta$-strong unfoldability
 embedding $j:M\to N$ with $j(\ell^*)(\kappa)=A$.
\end{enumerate}\label{ApproximateLD}

\Proof: Suppose that $\ell$ is a $\LDsunf_\kappa$ Laver Diamond function in $V$. Such a function has the Menas property, in the sense that
$j(\ell)(\kappa)$ can be made to have as large rank as desired, so we will assume that the lottery preparation has nontrivial forcing only at stages
$\gamma$ for which $\ell\image\gamma\of V_\gamma$. Let $\ell^*(\alpha)=x_{G_{\alpha+1}}$ when $\ell(\alpha)=\<x,y>$ for some $y$ and $x$ is a
$\P_{\alpha+1}$-name. For statement 1, suppose that $B\in P(\kappa)^{V[G]}$ and $A\in V_\theta[G][H]$. Let $\dot B$ be a name for $B$ of hereditary
size at most $\kappa$, and let $M$ be any $\kappa$-model in $V$ with $\dot B,\ell,\P\in M$. Let $\dot A$ be a $\P*\dot\Q$-name for $A$, and fix any
$\theta$-strong unfoldability embedding $j:M\to N$ with $j(\ell)(\kappa)=\<\dot A,\theta>$. The $\theta$ in the second coordinate ensures that the
next nontrivial stage of forcing is beyond $\theta$, and so the construction of the Main Theorem shows how to lift the embedding to
$j:M[G][H^*\intersect\Q^{M[G]}]\to N[j(G)][j(H^*\intersect\Q^{M[G]})]$, where $H^*$ is an isomorphic image of $H\of\Q$. Since $\dot B\in M$ it
follows that $B\in M[G]$. And since the construction has $j(G)_{\kappa+1}=G*H$ and $j(\ell)(\kappa)=\<\dot A,\theta>$, it follows by the definition
of $\ell^*$ that $j(\ell^*)(\kappa)=\dot A_{G*H}=A$. So the lifted embedding has all the features we desired for statement 1.

For statement 2, we are a bit more careful. Fix $B\in P(\kappa)^{V[G][H]}$ and $A\in V_\theta[G]$, with respective names $\dot A$ and $\dot B$ in
$V$. Let $\check{\dot A}$ be the result of viewing $\dot A$ as a $\P*\dot\Q$-name, rather than merely a $\P$-name. If there is no such
$\theta$-strong unfoldability embedding as in Statement 2, then there will be a condition $p\in\Q$ forcing over $V[G]$ that $\dot B$ cannot be placed
into a $\kappa$-model having a $\theta$-strong unfoldability embedding for which $j(\ell^*)(\kappa)=\check A$. Now fix as above any $\kappa$-model
$M$ containing $\dot B,\ell,\P$ and a $\theta$-strong unfoldability embedding $j:M\to N$ for which $j(\ell)(\kappa)=\<\check{\dot A},\theta>$. The
construction of the main theorem lifts this embedding to $j:M[G][H^*\intersect\Q^{M[G]}]\to N[j(G)][j(H^*\intersect\Q^{M[G]})]$, where $H^*$ is an
isomorphic image of $H\of\Q$, such that $p\in H^*$. Since $j(G)_{\kappa+1}=G*H$, it follows once again that $j(\ell^*)(\kappa)=(\check{\dot
A})_{G*H}=\dot A_G=A$. Let $B^*=\dot B_{G*H^*}$, which is determined of course by $G*(H^*\intersect\Q^{M[G]})$. Thus, we have constructed in
$V[G][H]=V[G][H^*]$ a $\kappa$-model $M[G][H^*\intersect\Q^{M[G]}]$ containing $B^*$ with an embedding for which $j(\ell^*)(\kappa)=A$. This
contradicts the fact that $p\in H^*$, since $p$ forced that there were no such model and embedding.\QED

It is not possible to combine statements 1 and 2 in our model, and have that for every $B\of\kappa$ in $V[G][H]$ and every $A\in V_\theta[G][H]$
there is a $\kappa$-model $M$ containing $B$ and $\ell^*$, with a $\theta$-strong unfoldability embedding $j:M\to N$ such that $j(\ell^*)(\kappa)=A$.
The reason this is not possible is that such a statement would exactly assert the Laver Diamond $\LDsunf_\kappa$ principle in $V[G][H]$, and this
implies $\Diamond_\kappa(\REG)$, which fails in $V[G][H]$. The point of Theorem \ref{ApproximateLD} is that one gets surprisingly close to
$\LDsunf_\kappa$ in $V[G][H]$. One can handle all $B$ and $A$, provided that at least one of the them is in $V[G]$. When $B\in V[G]$, then for any
$A\in V_\theta[G][H]$ we can ensure $B\in M[G][H^*\intersect\Q^{M[G]}]$ by placing a $\P$-name $\dot B$ into $M$. When $A\in V[G]$, then for
$B\of\kappa$ in $V[G][H]$ we can work with a $\P*\dot\Q$-name $\dot B$, and then interpret it by $G*H^*$, rather than by $G*H$, and conclude that
there can be no condition forcing that $\dot B$ has no embedding working with $A$. The point here is that if we would try to combine the two
arguments, and work in this case with an arbitrary $A\in V_\theta[G][H]$, with a name $\dot A$ in $V[G]$, then we would need to know that $\dot
A_{H^*}=\dot A_H$ in order to know that the lifted embedding, with $j(\ell^*)(\kappa)=A$ contradicts the fact that $p\in H^*$. In the actual
construction above, we have $A\in V[G]$ and so it has a name $\check\Adot$ that is not affected by $\pi$.

Statement 2 of Theorem \ref{ApproximateLD} of course includes the case when $A$ is an ordinal, and so in $V[G][H]$ we have an ordinal-anticipating
Laver diamond function for strong unfoldability.

\Corollary. For none of the principles $\LDwc_\kappa$, $\LDunf_\kappa$, $\LDind_\kappa$, $\LDsunf_\kappa$ does the existence of an
ordinal-anticipating Laver function necessarily imply the existence of a full set-anticipating Laver function. Indeed, the ordinal-anticipating Laver
functions do not even imply $\Diamond_\kappa(\REG)$.

\Proof: Our model shows that there can be an ordinal anticipating Laver function for strong unfoldability, the strongest of the four principles, but
no set-anticipating Laver function for weak compactness, the weakest notion. Since we have $\neg\Diamond_\kappa(\REG)$ in our model, the ordinal
anticipating Laver functions cannot provide $\Diamond_\kappa(\REG)$.\QED

\Corollary. If\/ $\LDsunf_\kappa$ holds in $V$, then in $V[G][H]$ there is a sequence\break $\<A_\alpha\st\alpha\in\REG\intersect\kappa>$ such that
\begin{enumerate}
 \item For any $A\of\kappa$ in $V[G]$, the set $\set{\alpha\st A\intersect\alpha=A_\alpha}$ is stationary in $V[G][H]$.
 \item For any $A\of\kappa$ in $V[G][H]$, the set $\set{\alpha\st A\intersect\alpha=A_\alpha}$ meets every club subset of $\kappa$ in $V[G]$.
\end{enumerate}

\Proof: Let $\ell^*$ be the function of Theorem \ref{ApproximateLD}, and let $A_\alpha=\ell^*(\alpha)$, if this is a subset of $\alpha$. If
$A\of\kappa$ is in $V[G]$, then for any club $C\of\kappa$ in $V[G][H]$ there is by part (2) of Theorem \ref{ApproximateLD} a $\kappa$-model $M$
containing $A$ and $C$ and an embedding $j:M\to N$ with $j(\ell^*)(\kappa)=A$. Since $\kappa\in j(C)$ and $j(A)\intersect\kappa=A=j(\ell^*)(\kappa)$,
it follows by reflection that $A\intersect\alpha=A_\alpha$ for many $\alpha\in C\intersect \REG$. In the case where $A\in V[G][H]$ and $C\in V[G]$,
then one argues similarly using part (1) of Theorem \ref{ApproximateLD}.\QED

These observations show that in any argument producing $\neg\Diamond_\kappa(\REG)$ for a strongly unfoldable cardinal $\kappa$, one should expect the
kind of circumlocutions that we went through in the Main Theorem with the unlucky case, rather than a pure lifting argument. To see why, suppose that
a pure lifting argument were able to succeed, meaning that every ground model embedding $j:M\to N$ lifts directly to $j:M[G][H\intersect\Q^M]\to
N[j(G)][j(H\intersect\Q^M)]$, where $j(G)=G*H*\Gtail$. Any desired $B\of\kappa$ in $V[G][H]$ could be included in $M[G][H\intersect\Q^M]$ just by
ensuring that $M$ was large enough. If $\ell$ were a ground model Laver function, then we could define $\ell^*(\alpha)$ to be
$\ell(\alpha)_{G_{\alpha+1}}$, and anticipate any $A\in V_\theta[G][H]$ by ensuring $j(\ell)(\kappa)=\dot A$ is a name for it. Thus, $\ell^*$ would
be a Laver function in $V[G][H]$, and we couldn't have $\neg\Diamond_\kappa(\REG)$ there, a contradiction.  It follows that the argument of the Main
Theorem must often be in the unlucky case, as this is the case that gives rise to the non-lifting circumlocutions.

Lastly, we close by mentioning that $\Diamond_\kappa(\REG)$ can be forced, while preserving the strong unfoldability of $\kappa$.

\Theorem.({\cite[Theorem 35]{Hamkins:LaverDiamond}}) If $\kappa$ is strongly unfoldable, then this is preserved to a forcing extension satisfying
$\LDsunf_\kappa$, which therefore also satisfies $\Diamond_\kappa(\REG)$.

\Proof: We quickly sketch the argument. The results of \cite{Hamkins2001:UnfoldableCardinals} show that the strong unfoldability of $\kappa$ is
preserved by fast function forcing $\F$, which adds a function $f\from\kappa\to\kappa$ such that every $(\theta+1)$-strong unfoldability embedding
$j:M\to N$ can be lifted to a $(\theta+1)$-strong unfoldability embedding $j:M[f]\to N[j(f)]$ in $V[f]$, and furthermore for any $\alpha<j(\kappa)$
such a lift can be found with $j(f)(\kappa)=\alpha$. Let $\vec a=\<a_\alpha\st\alpha<\kappa>$ be a fixed enumeration of $V_\kappa$, and let
$\ell(\gamma)=(a_{f(\gamma)})_{f\restrict\gamma}$, provided that $a_{f(\gamma)}$ is an $\F_\gamma$-name, where $\F_\gamma$ is the fast function
forcing at $\gamma$. This defines the Laver function $\ell\from\kappa\to V_\kappa$. Fix any $A\in V_\theta$, and select a name $\dot A$ so that
$A=\dot A_f$ and $\dot A\in V_\theta$. If $j:M\to N$ is a ground model $(\theta+1)$-strong unfoldability embedding, with $\vec a\in M$, then $\dot
A=j(\vec a)(\alpha)$ for some $\alpha$. We may lift the embedding to $j:M[f]\to N[j(f)]$ with $j(f)(\kappa)=\alpha$. It follows that
$j(\ell)(\kappa)=(j(\vec a)(\alpha))_{j(f)\restrict\kappa}=\dot A_f=A$, as desired. So $\LDsunf_\kappa$ holds in $V[f]$. This easily implies
$\Diamond_\kappa(\REG)$ by simply restricting $\ell$ to the values $\alpha$ for which $\ell(\alpha)\of \alpha$. For any $A\of\kappa$ and any club
$C\of\kappa$, there is a $\kappa$-model $M$ with $A,C,\ell\in M$ and $j:M\to N$ with $j(\ell)(\kappa)=A$. Since $\kappa\in j(C)\intersect\inacc$, it
follows by reflection that $\ell(\alpha)=A\intersect\alpha$ for some $\alpha\in C\intersect\inacc$. So $\Diamond_\kappa(\REG)$ holds.\QED

\bibliographystyle{alpha}
\bibliography{MathBiblio,HamkinsBiblio}

\end{document}